\documentclass{conm-p-l}
\newtheorem{theorem}{Theorem}

\theoremstyle{definition}

\theoremstyle{remark}

\def\range{\operatorname{range}}
\def\span{\operatorname{span}}
\def\rank{\operatorname{rank}}
\def\id{\operatorname{Id}}
\begin{document}
\title[The geometry of the skew-symmetric curvature operator]
{The Geometry of the Skew-Symmetric Curvature Operator in the
            Complex Setting}
\author[Gilkey]{Peter Gilkey}
\address{PG: Mathematics Department, University of Oregon, Eugene Or 97403 USA}
\email{gilkey@darkwing.uoregon.edu}
\author[Ivanova]{Raina Ivanova}
\address{RI: Dept. of Descriptive Geometry, University of Architecture,
Civil Engineering \& Geodesy, 1, Christo Smirnenski Blvd., 1421 Sofia,
Bulgaria}
\email{ivanovar@is.tsukuba.ac.jp and ivanovar@hopf.uoregon.edu}
\dedicatory{This paper is dedicated to Alfred Gray}
\begin{abstract} We shall construct almost complex algebraic curvature tensors for pseudo
Hermitian inner products whose skew-symmetric curvature operator has constant Jordan normal
form on the set of non-degenerate complex lines.
\end{abstract}
\maketitle
Let $(M,g)$ be a connected pseudo-Riemannian
manifold of signature
$(p,q)$ and dimension $m=p+q$. We are interested in local theory so we take
$M$ to be a small ball around a point $P$ of
$\mathbb{R}^m$ and take $g$ to be the germ of a metric defined near $P$. Let
${}^g\nabla$ be the Levi-Civita connection. The associated
curvature operator and curvature tensor are then given by:
\begin{eqnarray}\label{arefaa}
&&{}^gR(X,Y):={}^g\nabla_X{}^g\nabla_Y-{}^g\nabla_Y{}^g\nabla_X
             -{}^g\nabla_{[X,Y]},\text{ and}\\
&&{}^gR(X,Y,Z,W):=g({}^gR(X,Y)Z,W).\nonumber\end{eqnarray}
We have the
following symmetries:
\begin{eqnarray}  \label{arefab}
     &&{}^gR(X,Y,Z,W)=-{}^gR(Y,X,Z,W),\\
     &&{}^gR(X,Y,Z,W)={}^gR(Z,W,X,Y),\text{ and}\label{arefac}   \\
     &&{}^gR(X,Y,Z,W)+{}^gR(Y,Z,X,W)+{}^gR(Z,X,Y,W)=0.\label{arefad}
\end{eqnarray}
It is convenient to work in a purely algebraic setting. Let
$(\cdot,\cdot)$ be a symmetric bilinear form of signature
$(p,q)$ on a finite dimensional real vector space $V$. A $4$ tensor $R\in\otimes^4V$ is said
to be an  {\it algebraic curvature tensor} if
$R$ satisfies the symmetries given in equations (\ref{arefab})-(\ref{arefad}). A
pseudo-Riemannian manifold $(M,g)$ is said to be a {\it geometrical realization} of an
algebraic curvature tensor $R$ at a point $P$ of $M$ if there exists an isometry
$\Psi$ from the tangent space $T_PM$ to $M$ at $P$ to the vector space $V$ so that
$\Psi^*R={}^gR_P$ - i.e. we have:
\begin{eqnarray*}
&&(\Psi X,\Psi Y)=(X,Y)\text{ for all }X,Y\in T_PM\text{ and }\\
&&R(\Psi X,\Psi Y,\Psi Z,\Psi W)={}^gR(X,Y,Z,W)\text{ for all }X,Y,Z,W\in
T_PM.
\end{eqnarray*}
Since every algebraic curvature tensor is geometrically realizable by the germ of a
pseudo-Riemannian metric, algebraic curvature tensors have a central role in
differential geometry. To distinguish between the geometric and the algebraic contexts, we
shall use capital letters $X$, $Y$, $Z$, $W$ for tangent vectors and lower case letters $x$,
$y$, $z$, $w$ for elements of the vector space $V$.

If $\{x,y\}$ is an oriented basis for a non-degenerate
$2$ plane $\pi$, we define the {\it skew-symmetric curvature operator:}
$$R(\pi):=|(x,x)(y,y)-(x,y)^2|^{-1/2}R(x,y).$$
Let $J$ be a pseudo-Hermitian almost complex structure on $V$, i.e. $J$ is an
isometry with $J^2=-\id$. We use
$J$ to give a complex structure to $V$ by defining $\sqrt{-1}x:=Jx$. A
$2$ plane $\pi$ is said to be a {\it complex line} if $J\pi=\pi$. An algebraic curvature
tensor $R$ is said to be {\it almost complex} if
$JR(\pi)=R(\pi)J$ for every non-degenerate complex line $\pi$; this means that $R(\pi)$ is
complex linear.
\begin{theorem}\label{arefx} An algebraic curvature tensor $R$ is almost complex if and
only if $J^*R=R$ i.e. $R(x,y,z,w)=R(Jx,Jy,Jz,Jw)$ for all $(x,y,z,w)$ in $V$.
\end{theorem}

The algebraic curvature tensors satisfying the symmetry $J^*R=R$ form the class
${\mathcal{L}}_3$ in the classification of Gray \cite{refGray}; we also refer
Falcitelli, Farolina, and Salamon \cite{refFalcitelliFarinolaSalamon} for
related work. 

Since the curvature tensor encodes much of the geometry of the manifold, it can be a
very complicated object to study. There are natural operators associated to
$R$ which are  useful to be examined. One wants to know what are the geometric
consequences if the eigenvalues (or more generally the Jordan normal form) of such an
operator are constant on the natural domain of definition. We refer to the survey articles
\cite{refGilkey3, GilkeyIvanova1} for a more comprehensive introduction to the subject than
we can present here. In this paper, we shall study the geometry of the skew-symmetric
curvature operator for an almost complex algebraic curvature tensor.

Let $\phi$ be a
linear transformation of $V$. We use the inner product $(\cdot,\cdot)$ to define the adjoint
map
$\phi^*$ by the identity
$(\phi v,w)=(v,\phi^*w)$ for all $v$ and $w$. We say that $\phi$ is {\it self-adjoint} if
$\phi=\phi^*$ and that $\phi$ is {\it skew-adjoint} if $\phi=-\phi^*$. We shall need two
different families of algebraic curvature tensors.  If $\phi$ is self-adjoint, then we define:
\begin{eqnarray}
  &&R_\phi(x,y)z:=(\phi(y),z)\phi(x)-(\phi(x),z)\phi(y)\text{ and}\nonumber\\
  &&R_\phi(x,y,z,w):=(\phi(y),z)(\phi(x),w)-(\phi(x),z)(\phi(y),w).
     \label{arefae}\end{eqnarray}
If $\phi$ is skew-adjoint, then we define:
\begin{eqnarray}
  &&R_\phi(x,y)z:=(\phi y,z)\phi x-(\phi x,z)\phi y-2(\phi x,y)\phi z\text{ and}
   \nonumber\\
  &&R_\phi(x,y,z,w):=(\phi y,z)(\phi x,w)
  -(\phi x,z)(\phi y,w)-2(\phi x,y)(\phi z,w).
\label{arefaf}
 \end{eqnarray}

Let
$\mathbb{R}^{(p,q)}$ be Euclidean space with the standard metric of signature
$(p,q)$. Let  $S$ be the shape operator of a non-degenerate hypersurface in
$\mathbb{R}^{(p,q)}$. Then
${}^gR=\pm R_S$. Thus the tensor of equation (\ref{arefae}) arises in the study of
hypersurfaces in flat space. For example the algebraic curvature tensor
$\pm R_{\id}$ has constant sectional curvature
$\pm1$ and can be realized geometrically by the pseudo-spheres $S^\pm:=\{P:(P,P)=\pm1\}$.

In the Riemannian setting, $R_{\id}+R_J$ is the algebraic curvature
tensor of the Fubini-Study metric on $\mathbb{CP}^n$.
The following Theorem gives some of the algebraic
properties of the curvature tensors $R_\phi$ which we shall need.

\begin{theorem}\label{arefb}Assume that $\phi^*=\varepsilon\phi$ and that $J\phi=\varrho\phi
J$, where $\varepsilon=\pm1$ and $\varrho=\pm1$. Then
$R_\phi$ is an almost complex algebraic curvature tensor.
\end{theorem} 

Gray \cite{refGray} (Corollary 3.2) showed that if an algebraic
curvature tensor is geometrically realizable by a holomorphic Hermitian manifold (i.e. by a
manifold with an integrable almost complex structure), then $R$ satisfies the additional
symmetry:
\begin{eqnarray}
    &&R(x,y,z,w)+R(Jx,Jy,Jz,Jw)\nonumber\\
    &&\qquad=R(Jx,Jy,z,w)+R(Jx,y,Jz,w)+R(Jx,y,z,Jw)\nonumber\\
    &&\qquad+R(x,Jy,Jz,w)+R(x,Jy,z,Jw)+R(x,y,Jz,Jw).\label{arefBa}
\end{eqnarray}

\begin{theorem}\label{arefB} If $J\phi=\phi J$ and if $\phi^*=\pm\phi$, then $R_\phi$
satisfies equation {\rm(\ref{arefBa})}.\end{theorem}

\noindent{\bf Remark:} If $\rank(\phi)>2$, if $J\phi=-\phi J$ and if $\phi^*=\pm\phi$,
then $R_\phi$ does not satisfy equation (\ref{arefBa}); we omit details in the interests of
brevity.

\medbreak Two linear maps $T_1$ and $T_2$ of $V$ are said to be {\it Jordan equivalent} if
there exists an invertible linear map $\psi$ so that $T_1=\psi T_2\psi^{-1}$ or equivalently,
if there exists bases for $V$ so that the matrix representations of $T_1$ and $T_2$ with
respect to these two bases are equal. This means that the Jordan normal forms of $T_1$ and
$T_2$ are equal. Two such maps necessarily have the same eigenvalues but two maps with the
same eigenvalues need not be Jordan equivalent. In the positive definite setting $(p=0)$, the
Jordan normal form of a symmetric or of a skew-symmetric map is determined by the eigenvalue
structure. As the eigenvalue structure does not determine the Jordan normal form in the higher
signature setting ($p>0$), we work with the Jordan normal form rather than with the
eigenvalues.

An algebraic curvature tensor $R$ is said to be {\it Jordan IP} if the Jordan normal form of
the skew-symmetric curvature operator is constant on the Grassmannian of oriented
spacelike $2$ planes, is constant on the Grassmannian of oriented timelike $2$ planes, and is
constant on the Grassmannian of oriented mixed $2$ planes; the Jordan normal form is allowed
to vary with the type of the $2$ plane. One can show that if
$R$ is Jordan IP, then the rank $r$ is independent of the type of the plane. A
pseudo-Riemannian manifold
$(M,g)$ is said to be {\it Jordan IP} if
${}^gR$ is Jordan IP at every point $P$ of $M$; the Jordan normal forms are allowed to vary
with
$P$.

Suppose that $p=0$, i.e. that we are in the Riemannian setting. Ivanova and Stanilov
\cite{refIvanovaStanilov} began the study of the geometry of the skew-symmetric
curvature operator with constant eigenvalues (see also related work by Ivanova
\cite{refIvanova1, refIvanova2, refIvanova3}). Later Ivanov and Petrova
\cite{refIvanovPetrova} classified the Riemannian Jordan IP algebraic curvature
tensors and metrics in case $m=4$; for this reason the notation ``IP'' has
been used by some authors. Subsequently, the Riemannian Jordan IP algebraic curvature
tensors and metrics were classified in dimensions $m=5$, $m=6$, and
$m\ge9$ by Gilkey, Leahy, and Sadofsky \cite{refGilkeyLeahySadofsky} and in
dimension $m=8$ by Gilkey \cite{refGilkey2}; there are some partial results due
to Gilkey and Semmelman \cite{refGilkeySemmelman} if $m=7$ but the classification is
incomplete in this dimension. 

The first step in the study of Jordan IP algebraic curvature tensors is to control
the structure of the eigenspaces.  The following theorem was proved using topological methods
\cite{refGilkeyLeahySadofsky, refZhang}. It shows that the rank is $2$ in many cases. The
case $m=4$ is exceptional; there are Jordan IP algebraic curvature tensors of rank $4$ in
signatures $(4,0)$, $(2,2)$, and
$(0,4)$
\cite{refIvanovPetrova, refZhang}.
\begin{theorem}[Gilkey, Leahy, and Sadofsky; Zhang]\label{arefj}
Let $R$ be a Jordan IP algebraic curvature tensor of rank $r>0$ for a metric of signature
$(p,q)$.
\begin{enumerate}
\item Let $p\le1$. Let $q=5$, $q=6$, or $q\ge9$. Then $r=2$.
\item Let $p=2$. Let $q\ge10$. Let neither $q$ nor $q+2$ be
powers of $2$. Then $r=2$.
\end{enumerate}\end{theorem}

Theorem \ref{arefj} focuses attention on the Jordan IP algebraic curvature
tensors of rank
$2$. The following classification theorem follows from work of Gilkey
and Zhang \cite{refGilkeyZhang} extending previous work of Gilkey, Leahy, and Sadofsky
\cite{refGilkeyLeahySadofsky} from the Riemannian setting to arbitrary signature:

\begin{theorem}[Gilkey-Zhang]\label{arefk} Let $V$ have signature $(p,q)$, where $q\ge5$. 
Let $R$ be an algebraic curvature tensor on $V$. Then $R$ is a rank $2$ Jordan IP
algebraic curvature tensor if and only if there exists a non-zero constant $C$ so that
$R=CR_\phi$ where $\phi$ is self-adjoint and $\phi^2=\pm\id$.\end{theorem}

There is a corresponding classification of rank $2$ Jordan IP pseudo-Riemannian manifolds. We
refer to Gilkey and Zhang \cite{refGilkeyZhanga} for the proof of the following result:

\begin{theorem}[Gilkey-Zhang]\label{arefK} Let $(M,g)$ be a rank $2$ Jordan
IP pseudo-Rie\-mannian
manifold of signature $(p,q)$ for $q\ge5$. Then either $(M,g)$ has constant sectional
curvature or
$(M,g)$ is a warped product of an interval $I\subset\mathbb{R}$ with a
pseudo-Riemannian manifold $N$ of constant sectional curvature where
$$ds_M^2=\varepsilon dt^2+(\varepsilon Kt^2+At+B)ds^2_N\text{
and }\varepsilon=\pm1.$$\end{theorem}

Let $J$ be a pseudo-Hermitian almost complex structure on $V$. Let
$R$ be an almost complex algebraic curvature tensor. If the Jordan normal form of
$R(\cdot)$ is constant on the Grassmanian of non-degenerate complex lines, then $R$ is said
to be {\it almost complex Jordan IP}. 

If $p=0$, then the complex linear map $JR(\pi)$ is self-adjoint and thus
is diagonalizable with real eigenvalues. We refer to \cite{refGilkey4} for the proof
of the following result which controls the eigenvalue structure; the proof is topological in
nature.

\begin{theorem}[Gilkey]\label{arefm} Let $V$ have signature $(0,q)$. Let $R$ be
an almost complex Jordan IP algebraic curvature tensor on $V$. Let $\{\lambda_s,\mu_s\}$ be
the eigenvalues and multiplicities of $JR$, where we order the multiplicities $\mu_s$ so that
$\mu_0\ge...\ge\mu_\ell$. Suppose that
$\ell\ge1$. If $q\equiv2$ mod $4$, then $\ell=1$ and $\mu_1=1$. If $q\equiv0$ mod $4$, then
either $\ell=1$ and $\mu_1\le2$ or $\ell=2$ and $\mu_1=\mu_2=1$.
\end{theorem}

There are no general classification results for almost complex IP algebraic curvature
tensors and metrics. Kath \cite{refkath} has obtained some partial
results if
$(p,q)=(0,4)$ or $(p,q)=(2,2)$. Theorem \ref{arefn}
provides examples of almost complex Jordan IP algebraic curvature tensors. We adopt the
following notational conventions.
A linear transformation $\phi$ of $V$ is said to be {\it admissible} if
\goodbreak\begin{enumerate}
\item We have $\phi=\phi^*$ and $\phi J=\pm J\phi$ or $\phi=-\phi^*$ and $J\phi=-\phi J$.
\item We have $\phi^2=\id$, $\phi^2=-\id$, or $\phi^2=0$ and $\ker\phi=\range\phi$.
\end{enumerate}
A pair $\{\phi_1,\phi_2\}$ of linear transformations of $V$ is said to be {\it admissible}
if
\begin{enumerate}
\item Both $\phi_1$ and $\phi_2$ are admissible.
\item We have $\phi_1J=J\phi_1$, $\phi_2J=-J\phi_2$, and $\phi_1^*\phi_2+\phi_2^*\phi_1=0$.
\item If $\phi_1^2=\phi_2^2=0$ and if $\pi$ is any non-degenerate complex line, then we have
that
$\phi_1\pi\cap\phi_2\pi=\{0\}$.
\end{enumerate}

\begin{theorem}\label{arefn} Let $\phi$, $\phi_1$, and $\phi_2$ be linear transformations of
$V$, and let $c$, $c_0$, $c_1$, $c_2$, $c_3$ be arbitrary constants.
\begin{enumerate}
\item If $\phi$ is admissible, then $cR_\phi$ is an almost complex Jordan IP
algebraic curvature tensor.
\item If $\{\phi_1,\phi_2\}$ is admissible, then
$c_1R_{\phi_1}+c_2R_{\phi_2}$ is an almost complex Jordan IP algebraic curvature tensor.
\item The tensor $R:=c_0R_{\id}+c_1R_J$ is an almost complex Jordan IP algebraic curvature
tensor.
\item Let $\{i,j,k\}$ be a skew-adjoint quaternion structure on $V$. Then the tensor
$c_0R_{\id}+c_1R_i+c_2R_j+c_3R_k$ is an almost complex Jordan IP algebraic curvature
tensor relative to the almost complex Hermitian structure given by $J=i$.
\end{enumerate}
\end{theorem}

Let $J$ be the usual complex structure on
$\mathbb{R}^{(0,2s)}=\mathbb{C}^s$. Let $\{i,j,k\}$ be the usual skew-symmetric
quaternion structure on $\mathbb{R}^{(0,4s)}=\mathbb{C}^{2s}$. We can use Theorem
\ref{arefn} to construct examples realizing all the possibilities of Theorem \ref{arefm}. 

\begin{theorem}\label{arefp} Let $\lambda_s$ be distinct real numbers.
\begin{enumerate}
\item Let $R=c_0R_{\id}+c_1R_J$. We can choose $c_0$ and $c_1$ so
$JR(\cdot)$ has eigenvalues $\lambda_0$ and $\lambda_1$ with multiplicities $\mu_0=s-1$ and
$\mu_1=1$ on $\mathbb{C}^s$.
\item Let
$R=c_0R_{\id}+c_1R_i+c_2R_j+c_3R_k$. We can choose
$c_0$, $c_1$, $c_2$, and $c_3$ so
\begin{enumerate}
\item
$JR(\cdot)$ has eigenvalues $\lambda_0$ and $\lambda_1$ with multiplicities $\mu_0=2s-2$
and
$\mu_1=2$ on $\mathbb{C}^{2s}$ or
\item
$JR(\cdot)$ has eigenvalues $\lambda_0$, $\lambda_1$, and $\lambda_2$ with multiplicities
$\mu_0=2s-2$, $\mu_1=1$, and $\mu_2=1$ on $\mathbb{C}^{2s}$.
\end{enumerate}
\end{enumerate}
\end{theorem} 

\noindent The remainder of this paper is devoted to the proof of Theorems \ref{arefx},
\ref{arefb},
\ref{arefB}, \ref{arefn} and \ref{arefp}.

\begin{proof}[\bf Proof of Theorem \ref{arefx}] Let $R$ be an algebraic curvature tensor.
We follow an argument shown to us by Salamon. Note that $JR(\pi)=R(\pi)J$ for all
non-degenerate complex lines if and only if $JR(x,Jx)=R(x,Jx)J$ for all $x$ in $V$. We have
the following series of equivalences and implications:
\smallbreak\noindent\quad\phantom{$\Leftrightarrow$} a) $JR(x,Jx)=R(x,Jx)J$
$\forall x$,
\smallbreak\noindent\quad $\Leftrightarrow$ b) $(JR(x,Jx)z,w)-(R(x,Jx)Jz,w)=0$
$\forall x,z,w$,
\smallbreak\noindent\quad $\Leftrightarrow$ c) $R(x,Jx,z,Jw)+R(x,Jx,Jz,w)=0$
$\forall x,z,w$.
\smallbreak\noindent We polarize c), which is a quadratic identity in $x$, to obtain a
multilinear identity:
\smallbreak\noindent\quad $\Leftrightarrow$ d)
$R(y,Jx,z,Jw)+R(x,Jy,z,Jw)+R(y,Jx,Jz,w)+R(x,Jy,Jz,w)=0$
\smallbreak\noindent for all $x,y,z,w$. Replacing $(x,w)$ by $(Jx,Jw)$ and using
symmetry (\ref{arefab}) we obtain:
\smallbreak\noindent\quad $\Leftrightarrow$ e)
$-R(x,y,z,w)-R(Jx,Jy,z,w)+R(x,y,Jz,Jw)+R(Jx,Jy,Jz,Jw)=0$
\smallbreak\noindent for all $x,y,z,w$. Now symmetry (\ref{arefac}) implies:
\smallbreak\noindent\quad $\Leftrightarrow$ f)
$-R(z,w,x,y)-R(z,w,Jx,Jy)+R(Jz,Jw,x,y)+R(Jz,Jw,Jx,Jy)=0$
\smallbreak\noindent for all $x,y,z,w$. Interchange $(x,y)\leftrightarrow(z,w)$ to see:
\smallbreak\noindent\quad $\Leftrightarrow$ g)
$-R(x,y,z,w)-R(x,y,Jz,Jw)+R(Jx,Jy,z,w)+R(Jx,Jy,Jz,Jw)=0$
\smallbreak\noindent for all $x,y,z,w$. We add e) and g) to derive the identity:
\smallbreak\noindent\quad $\Rightarrow$ h) $-R(x,y,z,w)+R(Jx,Jy,Jz,Jw)=0$
    for all $x,y,z,w$
\smallbreak\noindent\quad $\Leftrightarrow$ i) $J^*R=R$. 
\smallbreak\noindent Conversely, we replace $(x,y)$ by
$(Jx,Jy)$ in h) and obtain:
\smallbreak\noindent\quad $\Leftrightarrow$ j) $-R(Jx,Jy,z,w)+R(x,y,Jz,Jw)=0$
for all $x,y,z,w$. \smallbreak\noindent 
Now adding h) and j), we derive e) which is equivalent to a).\end{proof}

\begin{proof}[\bf Proof of Theorem \ref{arefb}] Let $\phi^*=\pm\phi$. To prove that
$R_\phi$ is an algebraic curvature tensor, we must verify that
equations (\ref{arefab}), (\ref{arefac}), and (\ref{arefad}) are satisfied. Curvature
symmetry (\ref{arefab}) is immediate from the defining relations:
\begin{eqnarray*}
  &&R_\phi(x,y,z,w)=(\phi(y),z)(\phi(x),w)-(\phi(x),z)(\phi(y),w)\text{ if }\phi=\phi^*,
  \text{ and}\\
  &&R_\phi(x,y,z,w):=(\phi y,z)(\phi x,w)
  -(\phi x,z)(\phi y,w)-2(\phi x,y)(\phi z,w)\text{ if }\phi=-\phi^*.\end{eqnarray*}
We prove that $R_\phi$ satisfies equation (\ref{arefac}) by showing that the pairs $(x,y)$
and
$(z,w)$ play symmetric roles. Next, we verify that each of the terms comprising $R_\phi$ has
this symmetry by computing:
\begin{eqnarray*}
&&(\phi y,z)(\phi x,w)=(y,\phi z)(x,\phi w),\\
&&(\phi x,z)(\phi y,w)=(x,\phi z)(y,\phi w),\text{ and}\\
&&(\phi x,y)(\phi z,w)=(\phi z,w)(\phi x,y).\end{eqnarray*}
If
$\phi=\phi^*$, then we check that equation (\ref{arefad}) is satisfied by computing:
\medbreak\qquad $R_\phi (x,y,z,w)+ R_\phi (y,z,x,w) + R_\phi (z,x,y,w)$
\par\qquad\qquad$=(\phi y,z)(\phi x,w) - (\phi x,z)(\phi y,w)$
\par\qquad\qquad$ + \phantom{.}(\phi z,x)(\phi y,w)-(\phi y,x)(\phi z,w) $
\par\qquad\qquad$ +\phantom{.} (\phi x,y)(\phi z,w) - (\phi z,y)(\phi
x,w) =0$.
\medbreak\noindent If $\phi^*=-\phi$, then similarly we obtain:
\medbreak\qquad $R_\phi(x,y)z+R_\phi(y,z)x+R_\phi(z,x)y$
\par\qquad\qquad$=(\phi y,z)\phi x-(\phi x,z)\phi y-2(\phi x,y)\phi z$
\par\qquad\qquad$+\phantom{.}(\phi z,x)\phi y-(\phi y,x)\phi z-2(\phi y,z)\phi x$
\par\qquad\qquad$+\phantom{.}(\phi x,y)\phi z-(\phi z,y)\phi x-2(\phi z,x)\phi
y=0$.

\medbreak Let $\phi^*=\pm\phi$ and let $\phi J=\pm J\phi$. We have:
\begin{eqnarray*}
&&(\phi Jy,Jz)(\phi Jx,Jw)=(J\phi y,Jz)(J\phi x,Jw)
     =(\phi y,z)(\phi x,w),\\
&&(\phi Jx,Jz)(\phi Jy,Jw)=(J\phi x,Jz)(J\phi y,Jw)
     =(\phi y,z)(\phi x,w),\\
&&(\phi Jx,Jy)(\phi J z,Jw)=(J\phi x,Jy)(J\phi z,Jw)
     =(\phi x,y)(\phi z,w).\end{eqnarray*}
Since $R_\phi$ is a linear combination of these terms, we combine these three identities to
obtain
$R(Jx,Jy,Jz,Jw)=R(x,y,z,w)$. We use Theorem \ref{arefx} to see that $R$ is almost complex.
\end{proof}

\begin{proof}[\bf Proof of Theorem \ref{arefB}]
Let $\phi^*=\pm\phi$ and suppose $J\phi=\phi J$. By Theorem \ref{arefb}, $R_\phi$ is almost
complex. Thus by Theorem \ref{arefx} we see that:
\begin{eqnarray}
  &&R(x,y,z,w)=R(Jx,Jy,Jz,Jw),\ 
  R(Jx,Jy,z,w)=R(x,y,Jz,Jw),\nonumber\\ 
  &&R(Jx,y,Jz,w)=R(x,Jy,z,Jw),\ 
  R(Jx,y,z,Jw)=R(x,Jy,Jz,w).\label{arefxx}\end{eqnarray}
To prove Theorem \ref{arefB}, it suffices to show that equation (\ref{arefBa}) holds. Using
the relations of display (\ref{arefxx}), we represent equation (\ref{arefBa}) in the form:
\begin{equation}
R(x,y,z,w)=R(Jx,Jy,z,w)+R(Jx,y,Jz,w)+R(Jx,y,z,Jw).\label{arefxxz}\end{equation}
Equation (\ref{arefxxz}) is an additional symmetry. We complete the proof by showing that it
is satisfied by each of the terms comprising $R_\phi$:
\medbreak $(\phi Jy,z)(\phi Jx,w)+(\phi y,Jz)(\phi Jx,w)+(\phi y,z)(\phi Jx,Jw)
=(\phi y,z)(\phi x,w)$,
\smallbreak $(\phi Jx,Jz)(\phi y,w)+(\phi Jx,z)(\phi Jy,w)+(\phi Jx,z)(\phi y,Jw)
  =(\phi x,z)(\phi  y,w)$,
\smallbreak $(\phi Jx,Jy)(\phi z,w)+(\phi Jx,y)(\phi Jz,w)+(\phi Jx,y)(\phi z,Jw)
=(\phi x,y)(\phi z,w)$.\end{proof}

\begin{proof}[\bf Proof of Theorem \ref{arefn} (1)] Let $\phi$ be admissible. As
$\phi^*=\pm\phi$ and $J\phi=\pm\phi J$, $R_\phi$ is an almost complex algebraic curvature
tensor (see  Theorem \ref{arefb}). If
$\phi^*=-\phi$, then we also assume that
$J\phi=-\phi J$. Thus
$$(\phi x,Jx)=(-J\phi x,x)=(\phi Jx,x)=-(Jx,\phi x)=0$$
so the third term in equation (\ref{arefaf}) vanishes and we have:
\begin{equation}\label{DREFaa}
R_\phi(x,Jx)z=(\phi Jx,z)\phi x-(\phi x,z)\phi Jx\end{equation}
for both self-adjoint and skew-adjoint $\phi$.

Let $\pi:=\span\{x,Jx\}$ be a non-degenerate complex line where $x$ is either a unit spacelike
or a unit timelike vector. Then
$\pi$ is either spacelike or timelike, $\pi$ is not mixed. Because $\phi J=\pm J\phi$, the set
$$\phi\pi=\span\{\phi x,\phi Jx\}=\span\{\phi x,J\phi x\}$$ is a complex subspace of
$V$. To complete the proof, we must show the the Jordan normal form of $R_\phi(\pi)$ is
independent of $\pi$. We distinguish two cases.

First, suppose that $\phi^2=\varepsilon\id$ where $\varepsilon=\pm1$ and
$\phi^*=\varrho\phi$ where $\varrho=\pm1$. Then
$(\phi x,\phi y)=(\phi^*\phi x,y)=\varepsilon\varrho(x,y)=\pm(x,y)$ so either $\phi$ is
an isometry or $\phi$ is a para isometry. Thus $\phi\pi$ is either spacelike or timelike
since $\pi$ is either spacelike or timelike. We use equation (\ref{DREFaa}) to see that in
this case
$R_\phi(\pi)$ is a 90 degree rotation in the plane
$\phi\pi$ and $R_\phi(\pi)$ vanishes on $(\phi\pi)^\perp$. This implies that the Jordan
normal form of
$R_\phi(\pi)$ is independent of $\pi$.

Second, suppose that $\phi^2=0$. We then have
$$(\phi x,\phi
y)=(\phi^*\phi x,y)=\pm(\phi^2x,y)=0\text{ for any }x,y$$
 so the range of $\phi$ is totally
isotropic. Since
$\ker\phi=\range\phi$, the kernel of $\phi$ contains neither spacelike nor timelike vectors.
Because $\pi$ is either spacelike or timelike this means that
$\ker\phi\cap\pi=\{0\}$. Consequently, $\phi$ is an isomorphism from $\pi$ to $\phi\pi$. Thus
$\dim\phi\pi=2$ and
$\phi\pi=\span\{\phi x,J\phi x\}\subset\range\phi$ is a totally isotropic complex line. We use
equation (\ref{DREFaa}) to prove that in this case $R_\phi(\pi)^2=0$ and
$\rank(R_\phi(\pi))=2$. This implies that the Jordan normal form of $R_\phi(\pi)$ is
independent of
$\pi$. \end{proof}

\begin{proof}[\bf Proof of Theorem \ref{arefn} (2)] Let $R:=c_1R_{\phi_1}+c_2R_{\phi_2}$,
where
$\{\phi_1,\phi_2\}$ is an admissible pair. By assertion (1), $R_{\phi_1}$ and $R_{\phi_2}$ are
almost complex algebraic curvature tensors and hence $R$ is an almost complex algebraic
curvature tensor.

Let $\pi:=\span\{x,Jx\}$ be a non-degenerate complex line where $x$ is either a unit
spacelike or a unit timelike vector. We show that
$\{\phi_1x,\phi_1Jx,\phi_2x,\phi_2Jx\}$ is an orthogonal set by computing:
\begin{eqnarray*}(\phi_1x,\phi_1Jx)&=&(\phi_1x,J\phi_1x)=-(J\phi_1x,\phi_1x)
    =0,\\
   (\phi_1x,\phi_2x)&=&(\phi_2^*\phi_1x,x)=-(\phi_1^*\phi_2x,x)
    =-(\phi_2x,\phi_1x)=0,\\
(\phi_1x,\phi_2Jx)&=&-(J\phi_2^*\phi_1 x,x)
    =(\phi_2^*\phi_1Jx,x)=-(\phi_1^*\phi_2Jx,x)\\
&=&-(\phi_2Jx,\phi_1x)=0,\\
(\phi_1Jx,\phi_2x)&=&(\phi_2^*\phi_1Jx,x)=-(J\phi_2^*\phi_1x,x)
   =(\phi_2^*\phi_1x,Jx)\\
   &=&-(\phi_1^*\phi_2x,Jx)=-(\phi_2x,\phi_1Jx)=0,\\
    (\phi_1Jx,\phi_2Jx)&=&(\phi_2^*\phi_1Jx,Jx)
       =-(\phi_1^*\phi_2Jx,Jx)=-(\phi_2Jx,\phi_1Jx)=0,\\
    (\phi_2x,\phi_2Jx)&=&-(\phi_2x,J\phi_2x)=(J\phi_2x,\phi_2x)
    =0.\end{eqnarray*}
Thus the two complex lines $\phi_1\pi$ and $\phi_2\pi$ are
orthogonal complex lines. To complete the proof, we
must show that the Jordan normal form of $R(\pi)$ is independent of $\pi$. We distinguish
three cases:

First, suppose that $\phi_1^2=\pm\id$ and that $\phi_2^2=\pm\id$. In proving assertion (1), we showed
that $R_{\phi_s}(\pi)$ is a 90 degree rotation in the plane $\phi_s\pi$ and
$R_{\phi_s}(\pi)$ vanishes on $\phi_s\pi^\perp$ ($s=1,2$). Since $\phi_1\pi$ and
$\phi_2\pi$ are orthogonal non-degenerate complex lines, the Jordan normal form of $R(\pi)$
is independent of $\pi$.

Second, suppose that $\phi_1^2=\pm\id$ and that $\phi_2^2=0$ (the argument if $\phi_1^2=0$
and $\phi_2^2=\pm\id$ is similar and is therefore omitted). In proving assertion (1), we
showed that $R_{\phi_1}(\pi)$ is a 90 degree rotation in the plane $\phi_1\pi$, that
$\phi_1\pi$ is a timelike or spacelike complex line, and that
$R_{\phi_1}(\pi)$ vanishes on $\phi_1\pi^\perp$. We also showed that $\range R_{\phi_2}(\pi)$
is $\phi_2\pi$, that $\phi_2\pi$ is a totally isotropic complex line, and that
$R_{\phi_2}^2=0$.  Consequently, $\phi_1\pi\cap\phi_2\pi=\{0\}$ and the
Jordan normal form of
$R(\pi)$ is independent of $\pi$.

Third, suppose that $\phi_1^2=\phi_2^2=0$. In proving assertion (1), we showed that
$\range R_{\phi_s}(\pi)=\phi_s\pi$, that $\phi_s\pi$ is a totally isotropic
complex line, and that $R_{\phi_s}^2=0$ ($s=1,2)$. In this case, we have assumed
that $\phi_1\pi\cap\phi_2\pi=\{0\}$. It now follows that $\rank(R_\phi(\pi))=4$ and
that $R_\phi(\pi)^2=0$. These two properties determine the Jordan normal form of
$R_\phi(\pi)$ and one sees easily that it does not depend on
$\pi$.\end{proof}

\begin{proof}[\bf Proof of Theorem \ref{arefn} (3,4)] We consider the following two cases in
parallel: either let
$R:=c_0R_{\id}+c_1R_J$ or let
$R:=c_0R_{\id}+c_1R_i+c_2R_j+c_3R_k$, where $\{i,j,k\}$ is a skew-adjoint
quaternion structure and $J=i$. The first case is not always contained in the second one as
not all almost complex structures can be extended to quaternion structures. However the proof
follows the same lines.

The tensor $R$ is an almost complex algebraic curvature tensor (see Theorem \ref{arefb}). Let
$\pi:=\span\{x,ix\}$ be a non-degenerate complex line where $x$ is a unit spacelike vector;
the timelike case is similar. We shall complete the proof by showing that the Jordan
normal form of $R(\pi)$ is independent of $\pi$. Since
$\{i,j,k\}$ is a skew-symmetric quaternion structure, $\{x,ix,jx,kx\}$ is an orthonormal
set. It follows from equations (\ref{arefae}) and (\ref{arefaf}) that:
\begin{eqnarray}
&&R(\pi)x=-(c_0+3c_1)ix,\ R(\pi)ix=(c_0+3c_1)x,\nonumber\\
&&R(\pi)jx=-(2c_1-c_2-c_3)kx,\ R(\pi)kx=(2c_1-c_2-c_3)jx,\label{arefz}\\
&&R(\pi)y=-(2c_1)iy\text{ if }y\perp\span\{x,ix,jx,kx\}.\nonumber\end{eqnarray}
Let $V_1:=\span\{x,ix\}$,
$V_2:=\span\{jx,kx\}$, and
$V_3:=(V_1+V_2)^\perp$. We have an orthogonal direct sum decomposition
$V=V_1\oplus V_2\oplus V_3$. These subspaces
are preserved by the almost complex structure $i$ and by $R(\pi)$. Let $R_s$ be the
restriction of
$R(\pi)$ to
$V_s$ for $1\le s\le 3$. It follows from equation (\ref{arefz}) that
$R_1=-(c_0+3c_1)i$, $ R_2=-(2c_1-c_2-c_3)i$, and $R_3=-(2c_1)i$. These are multiples
of positive rotations through 90 degrees. If $x$ is timelike, then the rotations are negative.
In this way we obtain the
Jordan normal form of $R(\pi)$ and verify it does not depend on $\pi$.
\end{proof}

\begin{proof}[\bf Proof of Theorem \ref{arefp}] Let $R=c_0R_0+c_1R_i+c_2R_j+c_3R_k$ on
$\mathbb{C}^{2s}$ and let $J=i$.
In proving Theorem \ref{arefn}, we determined the Jordan normal form of
$R(\pi)$. It now follows that $(c_0+3c_1)$ is an eigenvalue of multiplicity $1$, that
$(2c_1-c_2-c_3)$ is an eigenvalue of multiplicity $1$, and that $2c_1$ is an eigenvalue of
multiplicity $2s-2$ for the self-adjoint operator $JR(\pi)$. We choose the constants $c_s$
for $s=0,1,2,3$ so that $\lambda_1=c_0+3c_1$, $\lambda_2=2c_1-c_2-c_3$, and $\lambda_0=2c_1$
to complete the proof of the second assertion; the first follows similarly.
\end{proof}

\medbreak\noindent{\bf Acknowledgements.} It is a pleasure for the first author to acknowledge
helpful conversations with L. Vanhecke and S. Salamon concerning curvature
decompositions. He  is also grateful to A. Gray
who first introduced him to the subject of almost Hermitian pseudo-Riemannian geometry and
curvature decompositions.

The research of the first author was supported in part by the NSF
(USA) and the MPI (Leipzig, Germany). The research of the second author was supported in part
by the JSPS Post Doctoral Fellowship Program (Japan) and the MPI (Leipzig, Germany). 

\bibliographystyle{amsalpha}

\end{document}